\magnification 1200
        \def\R{{\rm I\kern-0.2em R\kern0.2em \kern-0.2em}}
        \def\N{{\rm I\kern-0.2em N\kern0.2em \kern-0.2em}}
        \def\P{{\rm I\kern-0.2em P\kern0.2em \kern-0.2em}}
        \def\B{{\rm I\kern-0.2em B\kern0.2em \kern-0.2em}}
        \def\Z{{\rm I\kern-0.2em Z\kern0.2em \kern-0.2em}}
        \def\C{{\bf \rm C}\kern-.4em {\vrule height1.4ex width.08em depth-.04ex}\;}
        \def\B{{\bf \rm B}\kern-.4em {\vrule height1.4ex width.08em depth-.04ex}\;}
        
        \def\D{{\Delta}}

        \def\z{{\zeta}}
        \def\cC{{\cal C}}

        \def\L{{\cal L}}
        \def\cA{{\cal A}}
        \def\cB{{\cal B}}
        \def\cC{{\cal C}}
        \def\cD{{\cal D}}
        \def\cS{{\cal S}}
        \def\cS{{\cal S}}
        
        \def\cV{{\cal V}}

        \def\cU{{\cal U}}
        \def\cT{{\cal T}}

        \def\tau{{\sigma}}

        \
        \vskip 20mm
        \centerline {\bf A COMPLETE COMPLEX HYPERSURFACE IN THE BALL OF $\C^N$}
        \vskip 4mm
        \centerline{Josip Globevnik}
        \vskip 4mm
        \bf Abstract\ \ \rm In 1977 P.\ Yang asked whether there exist complete immersed complex 
        submanifolds $\varphi \colon\ M^k\rightarrow \C^N$ with bounded image. A positive answer is known for
        holomorphic curves $(k=1)$ and partial answers are known for the case when $k>1$.
        The principal result of the present paper is a construction of a holomorphic function 
        on the open unit ball 
        $\B_N$ of $\C^N$ whose real part is unbounded on every path in $\B_N$ 
        of finite length that ends on $b\B_N$. A
        consequence is the existence of a complete, closed complex hypersurface in $\B_N$. This gives a positive 
        answer to 
        Yang's question in all dimensions $k,\ N,\ 1\leq k<N$, by providing properly embedded complete
        complex manifolds.
        \vskip 4mm
        \bf 1. Introduction and the main result \rm
        \vskip 2mm
        Denote by $\D$ the open unit disc in $\C$ and by $\B_N$ the open unit ball in $\C^N, \ N\geq 2$.
        
        In 1977 P.\ Yang asked whether there exist complete immersed complex submanifolds 
        $\varphi\colon M^k\rightarrow \C^N$ with bounded image [Y1, Y2]. The first answer 
        was obtained by P.\ Jones [J] who constructed 
        a bounded complete immersion $\varphi\colon\D\rightarrow \C^2$ and a complete
        proper holomorphic embedding 
        $\varphi\colon\D\rightarrow \B_4$. Since 
        then there has been a series of 
        results on bounded complete holomorphic curves ($k=1$) \it immersed \rm
        in $\C^2 $ [MUY, AL1, AF] the most recent being that every 
        bordered Riemann surface admits a complete proper holomorphic 
        immersion to $\B_2$ and a complete proper holomorphic embedding to $\B_3$ [AF]. 
        The more difficult complete \it embedding \rm problem for $k= 1$ and $N=2$ 
        has been solved only very  
        recently by A.\ Alarc\'{o}n and F.\ J.\ L\'{o}pez [AL2] who proved that every convex domain in $\C^2$ 
        contains a complete, properly embedded complex curve.
        
        In the present paper we are interested primarily in
        the higher dimensional case $(k>1)$ where there are partial answers which are 
        easy consequences of the results for complete curves. For instance, it is known that
        for any $k\in\N$ there are complete bounded embedded complex $k$-dimensional submanifolds of 
        $\C^{2k}$ and it is an open question whether, in this case,  $N=2k$ is 
        the minimal possible dimension [AL2]. 
        In the present paper we consider the case where 
        $\varphi $ is a proper holomorphic embedding. In this
        case $\varphi (M^k)$ is a closed submanifold. 
        We restate the definition of completeness for this case:
        \vskip 2mm
        \noindent \bf DEFINITION 1.1 \it A closed complex submanifold $M$ of $\B_N$ is complete if 
        every path  $p\colon [0,1)\rightarrow M$ such that $|p(t)| \rightarrow 1$ as 
        $t\rightarrow 1$ has infinite length. \rm
        \vskip 2mm
        \noindent Note that this coincides with the standard definition of completeness 
        since the paths
        
        \noindent $p\colon [0,1)\rightarrow M$ such that  $|p(t)| \rightarrow 1$ 
        as $t\rightarrow 1$ are precisely the paths that leave every compact subset 
        of $M$ as $t\rightarrow 1$ . 
        \vskip 1mm
        
        \noindent Here is our main result:
        \vskip 2mm
        \noindent \bf THEOREM 1.1\ \it Let $N\geq 2$. There is a holomorphic function $f$ 
        on $\B_N$ such that $\Re f$ is unbounded on every 
        path of finite length that ends on $b\B_N$.
        \vskip 2mm \rm 
        \noindent So our function $f$ has the property that 
        if $p\colon [0,1]\rightarrow \overline{\B_N}$ is a 
        path of finite length such that 
        $|p(t)|<1 \ (0\leq t <1)$ and $| p(1)| = 1$ then $t\rightarrow \Re \bigl(f(p(t)\bigr)$ 
        is unbounded on $[0,1)$.
        \vskip 1mm
         The following corollary answers the question of Yang in all
         dimensions $k$ and $N$ by providing properly embedded complete complex manifolds.
        \vskip 2mm
        \noindent
        \bf COROLLARY 1.2\ \it For each $k, N,\ 1\leq k<N,$ there is a complete, 
        closed, k-dimensional complex submanifold of  $\B_N$. \rm
        \vskip 1mm
        \noindent \bf Proof. \rm We first prove the corollary 
        for $k=N-1$ (that is, we first prove the existence of
        the hypersurface, mentioned in the title).
        Let $f$ be the function given by Theorem 1.1. By Sard's
        theorem one can choose $c\in\C$ such that the level set 
        $M=\{ z\in \B_N\colon f(z) = c\} $ is a closed submanifold of $\B_N$. Let
        $p\colon[0,1)\rightarrow M$ 
        be a path such that $p(t)\rightarrow b\B_N$ as $t\rightarrow 1$. Assume 
        that $p$ has finite length. Then there is a point $w$ on $b\B_N$ 
        such that $\lim_{t\rightarrow 1}p(t)= w$. By the properties of $f$,\ $\Re f$ 
        is unbounded on $p([0,1))$. On the other hand, $f((p(t))=c \ (0\leq t<1)$,
        a contradiction. So $p$  must have infinite length. This proves that $M$ is complete and
        so completes the proof of the corollary for $k=N-1$. Assume now that $1\leq k\leq N-2$. By
        the first part of the proof there is a complete, closed, $k-$dimensional 
        complex submanifold $M$ of
        $\B_{k+1}\subset \B_N$. Clearly $M$ is a complete, closed $k-$dimensional manifold of $\B_N$. 
        This completes the proof.
        \vskip 2mm
        \noindent\bf Remark\ \ \rm If we want to have a connected, complete closed complex 
        submanifold of $\B_N$ then we simply take a connected component of $M$ as above. Note also that
        the same function $f$ gives many complete closed complex manifolds of $\B_N$ since, by Sard's theorem, 
        one can use the same reasoning for almost every $c$ in the range of $f$.
        \vskip 4mm
        \bf 2.\ Outline of the proof of Theorem 1.1 \rm
        \vskip 2mm
        Let $M\in\N$. For $x\in\R^M\setminus\{ 0\}$ and $\alpha\in\R$, write
        $$ H(x,\alpha )=\{ y\in\R^M\colon <y|x>=\alpha\}, \ \ K(x,\alpha )=\{ y\in\R^M\colon <y|x>\leq\alpha\} .
        $$
        Assume that $x_i\in\R^M\setminus\{ 0\}\ \ (1\leq i\leq n)$ \ and that 
        $$
        P = \bigcap _{i=1}^n K(x_i,1)
        \eqno (2.1)
        $$
        is a bounded set. Then $P$ is a \it convex polytope,\rm \ that is, the convex hull 
        of a finite set.
        So $P$ is a compact convex set that contains the origin in its interior. A convex 
        subset $F$ of $P$ is called a 
        \it face \rm of $P$ if any closed segment with endpoints in $P$ whose relative 
        interior meets $F$ is contained in $F$. A
        \it k-face \rm is a face $F$ with $\hbox{dim} F = k$, that is, the affine hull
        of $F$ is k-dimensional. A face of dimension 
        $M-1$ is called a \it facet \rm of $P$. Let $P$ be a convex polytope such that 
        the representation (2.1) is 
        irreducible, that is,
        $$
        P\not= \bigcap_{i=1, i\not=k}^n K(x_i, 1)\hbox{\ \ for each\ \ }k,\ 1\leq k\leq n .
        $$
        Then 
        $$
        bP = \bigcup_{i=1}^n H(x_i, 1)\cap P
        $$ and the sets $F_i = H(x_i, 1)\cap P,\ 1\leq i\leq n$, 
        are precisely the facets of $P$. See [B] for the details.
        
        Given a convex set $G$, denote by $\hbox{ri} (G) $ the relative interior of $G$ in the 
        affine hull of $G$.
        What remains of the boundary of a convex polytope $P$ 
        after we have removed relative interiors of all facets $F_i,\ 1\leq i\leq n$, we call the 
        \it skeleton \rm of $P$ (or more precisely, the $(M-2)$-skeleton of $P$, 
        the union of all $(M-2)$-dimensional faces of $P$) and denote by 
        $\hbox{skel} (P)$. Thus
        $$
        \hbox{skel} (P)= \bigcup _{i=1}^n\bigl[ F_i\setminus \hbox{ri} (F_i)\bigr] .
        $$
        
        To prove Theorem 1.1 we first prove 
        \vskip 2mm
        \noindent\bf THEOREM 2.1 \it Let $\B $ be the open unit ball of $\R^M,\ M\geq 3$.
        There is a sequence of convex polytopes $P_n,\ n\in\N$, 
        such that
        $$ P_1\subset \hbox{\rm Int} P_2\subset P_2\subset \hbox{\rm Int}  P_3\subset\cdots \subset \B,
        \ \ \ \bigcup_{j=1}^\infty P_j = \B ,
        $$
        such that if $w_j\in\hbox{\rm skel}(P_j)\ \ (j\in\N)$ then 
        $$
        \sum_{j=1}^\infty
        |w_{j+1}-w_j| = \infty \eqno (2.2)
        $$  
        that is, the series in (2.2) diverges. 
        \vskip 2mm
        \noindent \rm In the proof of Theorem 1.1 we shall use the following 
        \vskip 2mm
        \noindent\bf COROLLARY 2.2\ \it Let $P_n,\ n\in\N$ be the sequence of convex polytopes from
        Theorem 2.1. Let $\theta_n$
        be a decreasing sequence of positive numbers such that
        $\sum_{n=1}^\infty \theta_n<\infty$. For
        each $n\in \N$, let $\cU_n\subset bP_n$ be the $\theta_n$-neighbourhood 
        of $\hbox{\rm skel}(P_n)$ in $bP_n$, that is, $\cU_n=\{ w\in bP_n\colon
        \ \hbox{\rm dist}(w, \hbox{\rm skel}(P_n))<\theta_n\}$.   
        Let $p\colon \ [0,1)\rightarrow \B$ be a path such that 
        $|p(t)|\rightarrow 1$ as $t\rightarrow 1$, 
        and such that for all sufficiently large $n\in\N$,\ \ 
        $p([0,1))$ meets $bP_n$ only at $\cU_n$.  Then
        $p$ has infinite length. \rm
        \vskip 2mm
        Once we have proved Corollary 2.2 we prove Theorem 1.1 as follows. Let $\B_N$ be 
        the open unit ball of $\C^N,\ N\geq 2$. 
        Let $P_n, \ n \in\N$, be a sequence of convex polytopes as in Theorem 2.1 with
        $M=2N$, and let $\cU_n, \ n\in\N$, be as in Corollary 2.2.  Given $\varepsilon _n>0$ 
        and $L_n<\infty$ we use an idea from [GS] to construct a function $f_n$, holomorphic 
        on $\B_N$, such that 
        $|f_n|<\varepsilon_n $ on $P_{n-1}$ and such that
        $\Re f_n >L_n$ on $ bP_n\setminus \cU_n$. By choosing 
        $L_n$ and $\varepsilon _n$ inductively in the right way we then see that
        $f=\sum_{n=1}^\infty f_n$ has all the required properties. 
        \vskip 4mm
        \bf 3.\ Beginning of the proof of Theorem 2.1\rm
        \vskip 2mm
        Let $w_n$ be a sequence in $\B$ such that $|w_n|\rightarrow 1$ as $n\rightarrow \infty$. If $w_n$ 
        does not converge then (2.2) holds so to prove Theorem 2.1 
        it is enough to consider only the convergent sequences $w_n$. 
         
        First, we try to explain the idea of the most important part of the proof. 
        Suppose for a moment that we have a sequence $P_n$ of convex polytopes with 
        the desired properties and that there is an increasing sequence $R_n$ of positive numbers 
        converging to $1$ such that
        $$
        bP_n\subset R_n\overline \B \setminus R_{n-1}\overline \B \ \ (n\in\N).
        $$
        Let $W=U\times (1-\nu, 1+\nu)$ be a small open
        neighbourhood of $z=(0,0,\cdots 0,1)$ in $\R^M$ where $U$ 
        is a small open ball in $\R^{M-1}$ centered at the origin and $\nu >0$ is small. Assume that $
        U\times\{ 1-\nu\} \subset R_0\B$.
        
        Let $\pi$ be the orthogonal projection onto $\R^{M-1}$, so 
        $\pi (x_1,\cdots, x_M) = (x_1, \cdots , x_{M-1})$. For each $n$, consider
        $C_n$, the part of $bP_n\cap W$ consisting of the facets of $P_n$ contained in $W$. The projection 
        $\pi$ 
        is one to one on $C_n$ and 
        for each of these facets its image under $\pi$ is a convex polytope in $U$ that is a cell of a 
        partition of $\pi (C_n)$ into convex polytopes. Call this partition $\L_n$ and notice that 
        as $n\rightarrow \infty$,  
        $\pi (C_n)$ tends to $U$. If we remove from each cell of $\L_n$ its 
        relative interior 
        then we get what we call the skeleton of $\L_n$ and denote by $\hbox{skel} (\L_n)$. Clearly $\pi
        (\hbox{skel} (P_n)\cap C_n ) = \hbox{skel} (\L_n)$. 
        Since, by our assumption at the moment, every sequence $w_n$ contained in $W$ which meets $\hbox{skel}(P_n)$ for 
        all sufficiently large $n$ must satisfy (2.2), looking at $z_n=\pi (w_n)$ we conclude that 
        every sequence $z_n\in U$ such that $z_n\in\hbox{skel}(\L_n)$  for all sufficiently 
        large  $n$ must satisfy $\sum_{n=1}^\infty |z_{n+1}-z_n|=\infty $.
        The idea now is to reverse the direction of reasoning. Let $R_0$ be so close to $1$ that
        $U\times\{1-\nu\}\subset R_0\B$. In a typical induction step of constructing 
        our polytopes the data will be a partition $\L$ of $\R^{M-1}$ into convex 
        polytopes and $\rho$ and $r$, 
        $R_0<\rho<r<1.$ Denote by $\cC$ the union of those cells of the partition $\L$ that
        are contained 
        in $U$ and let $\cV$ be the set of their vertices. We will "lift" 
        $\cV$ to $b(r\B)$ by putting $V=
        (\pi|W\cap b(r\B))^{-1}(\cV)$. We want $V$ to be the set of vertices of a convex 
        polyhedral surface 
        $C$ such that $\pi (C)=\cC$ and such that $\pi$ maps the facets of $C$ 
        precisely onto the cells of $\cC$. 
        We will do this in such a way that $C$ stays out of $\rho\overline \B$ - 
        for this, the cells of $C$, and 
        consequently the cells of $\cC$ will have to be sufficiently small, of
        size proportional to $\sqrt{r-\rho}$. 
        Then we will construct a convex polytope $P$ such that $C$ will be a part
        of its boundary $bP$ 
        and such that $\rho\overline\B\subset \hbox{Int} P\subset P\subset r\overline\B$.
        
        There is a potential problem already at the first step. Namely,  the points of $V$ need 
        not be the vertices of a \it convex\rm\  surface $C$. For this to happen we will need 
        two things:
        $\L$ will have to be a true Delaunay partition of $\R^{M-1} $ and the ball $U$ in 
        the definition of $W$ will have to be sufficiently small so that the part of $b(r\B) $ 
        contained in $W$ will be sufficiently flat. 
        \vskip 4mm
        \bf 4. A Delaunay tessellation of $\R^{M-1}$ \rm
        \vskip 2mm
        Perturb the canonical orthonormal basis in $\R^{M-1}$ a little to get an 
        $(M-1)$-tuple of vectors $e_1, e_2, \cdots , e_{M-1}$ in general position so that the lattice
        $$
        \Lambda = \bigl\{ \sum_{i=1}^{M-1}n_ie_i \ \colon \ n_i\in Z, \ 1\leq i\leq M-1\bigr\} 
        \eqno (4.1)
        $$
        will be generic, and, in particular, no more than $M$ points of $\Lambda$ 
        will lie on the same sphere. 
        
        For each point $x\in\Lambda$ there is the\  
        \it Voronei cell \rm\  $V(x)$ consisting of those points of
        $\R^{M-1}$ that are at least as close to $x$ as to any other $y\in\Lambda$, so
        $$
        V(x))=\{ y\in\R^{M-1}\colon\ \hbox{dist} (y,x)\leq \hbox{dist} (y,z)\hbox{\ for all\ }z\in\Lambda
        \} .
        $$In our case it is easy to see how to get $V(0)$.  Consider the finite set 
        $E=\{ \sum_{j=1}^{M-1}n_i e_i: -1\leq n_i\leq 1, 1 \leq i\leq M-1\} $ and for 
        each $x\in E\setminus \{ 0\}$, look at $K(x,|x|^2/2)$, that is, at the 
        halfspace which contains the origin and is 
        bounded by the hyperplane passing through $x/2$ which is perpendicular to $x$. Then
        $$
        V(0)=\bigcap_{x\in E\setminus\{ 0\}} K(x, |x|^2/2).
        $$
        This is a convex polytope. It is known that the Voronei cells form a tessellation of 
        $R^{M-1}$ and in our case they are all congruent, 
        of the form $V(0)+x, \ x\in\Lambda$ [CS]. 
        
        There is a \it Delaunay cell\rm\  for each point that is a vertex of a Voronei cell.
        It is the convex polytope that is the 
        convex hull of the points in $\Lambda $ closest to that point - these points are all on 
        a sphere centered at this point. In our case, when there are no more than $M$ 
        points of $\Lambda$ on a sphere, Delaunay cells are $(M-1)$-simplices. 
        Delaunay cells form a tessellation of $R^{M-1}$ [CS]. It is a \it true \rm 
         Delaunay tessellation , that is, for each cell, the circumsphere of each cell $S$ contains 
        no other points of $\Lambda$ than the vertices of $S$. We shall denote by $\cD (\Lambda )$ 
        the family of all simplices - cells of the Delaunay tessellation for the
        lattice $\Lambda $. 
        
        By periodicity there are finitely many simplices $S_1,\cdots, S_\ell$ such that
        every other simplex of $\cD (\Lambda )$ is of the form $S_i + w$ 
        where $w\in\Lambda$ and
        $1\leq i\leq \ell$. It is then clear by periodicity that there is an $\eta>0$ such that 
        for every simplex $S\in\cD (\Lambda )$  
        in $\eta$-neighbourhood 
        of the closed ball bounded by the circumsphere of $S$ there are no other points of  $\Lambda$ 
        than the vertices of $S$.

        We shall typically replace the lattice $\Lambda $ by the 
        lattice $\Lambda + q =\{ x+q\colon\ x\in\Lambda\} $ where $q\in\R^{M-1}$,  
        or, more generally, by the lattice $\tau(\Lambda + q)$ 
        where $\tau>0$ is small. Again, we shall denote 
        by $\cD(\tau (\Lambda + q))$ the family of all simplices - cells of 
        the Delaunay tessellation for $\tau (\Lambda+q)$. These are the simplices of the form $\tau (S+q)$ where 
        $S\in \cD (\Lambda)$. Passing from $\Lambda$ to $\tau(\Lambda +q)$ everything in the reasoning
        will change proportionally.
        In particular, for every simplex 
        $S\in \cD(\tau(\Lambda + q)) $ in $(\tau\eta)$-neighbourhood of the closed ball bounded by the 
        circumsphere of $S$ there will be no other points of $\tau(\Lambda + q)$ than the vertices 
        of $S$.  We shall also need the notion of the \it skeleton \rm of the 
        Delaunay tessellation for $\tau(\Lambda + q)$. This is what remains after we 
        remove the interiors of all 
        $ S\in\cD (\tau(\Lambda + q))$, hence
        $$
        \hbox{skel}\bigl(\cD (\tau(\Lambda +q))\bigr) = \bigcup_{S\in\cD (\tau(\Lambda + q))}
        \bigl[S\setminus\hbox{Int}S\bigr] = \R^{M-1}\setminus \bigl[\bigcup _{S\in\cD 
        (\tau(\Lambda + q))}\hbox{Int}S\bigr].
        $$

        The author is gratefull to John M.Sullivan who suggested the use of a generic lattice
        for our purpose here. 
        \vskip 4mm
        \bf 5. Lifting the lattice from $\R ^{M-1}$ to the sphere \rm
        \vskip 2mm
        Let $z$, $W=U\times (1-\nu, 1+\nu )$ and $\pi$ be as in Section 3. Let $\Lambda\subset\R^{M-1}$ 
        be as in (4.1). 
        
        Fix $R_0,\ 0< R_0<1$, so large that $U\times \{ 1-\nu\}\subset R_0\B$ and assume that 
        $R_0<\rho<r<1$. The part of the sphere $b(r\B ) $ in $W$ can now be written as a graph 
        of a real analytic 
        function, call it $\psi_r$, so
        $$
        b(r\B)\cap W = \{ (x,\psi_r (x))\colon\ x\in U\}
        $$
        where
        $$
        \psi_r (x)=\psi_r (x_1,\cdots ,x_{M-1}) = \biggl( r^2-\sum_{j=1}^{M-1}x_j^2\biggr)^{1/2} .
        \eqno (5.1)
        $$
        Note that $(\hbox{grad}\ \psi_r)(0)=0,\ R_0<r<1$.
        
        The map $\pi $ maps $W\cap b(r\B)$ in a one to one
        way onto $U$. We shall "lift" 
        $(\tau\Lambda)$ from $U $ to $b(r\B)\cap W$ by the inverse of 
        this map, that is, by the map $x\mapsto (x,\psi_r(x))$.  We want to get a 
        convex polyhedral surface $C$
         with vertices $w=(v,\psi_r (v))$ where $v$ are the vertices of those 
         cells of the Delaunay tessellation for 
         $\tau\Lambda$ which are contained in $U$ and we want that $\pi$ maps
         the facets of the surface $C$ precisely onto the Delaunay
          cells of $\tau\Lambda$ contained in $U$. Let us describe the 
          conditions for this 
          to happen. Let $S$ be a simplex of the Delaunay 
        tessellation for $\tau\Lambda$. Let $v_1,\cdots , v_M$ be the vertices of 
        $S$. We want that the simplex with 
        vertices $w_j = (v_j, \psi (v_j)),\ 1\leq j\leq M$, is a facet of a
        convex poyhedral surface.
        For this to happen, 
        all other points $w=(v,\psi_r (v)), \ v\in\tau\Lambda \cap U,\  v\not= v_1,\cdots, v_M,$ 
        must lie in 
        the open halfspace bounded by the hyperplane $\Pi$
        through $w_j,\ 1\leq j\leq M$, which contains the origin, that is, 
        they must lie on $b(r\B)$ outside the "small"
        sphere $\Gamma 
        = \Pi\cap b(r\B)$. Since $\pi|W\cap b(r\B)$ is one to one, 
        this happens if and only if the points $v\in\tau\Lambda$ which are the 
        vertices of the Delaunay cells of $\tau\Lambda$ contained in $U$ and are different 
        from $v_1,\cdots , v_M$, are outside the projection $\pi (\Gamma)$, an ellipsoid in $\R^{M-1}$.
        
        As we shall see, this will happen for all such simplices $S$ if the ball $U\subset \R^{M-1}$ 
        centered at the origin 
        will be small 
        enough so that the the gradient of $\psi_r $ and thus the Lipschitz constant of $\psi_r $ 
        will be small enough on 
        $U$. The choice of $U$ will depend only on $\eta$ from Section 4 and 
        the same reasoning will work for any $\tau >0$.
        \vskip 2mm
        \noindent\bf LEMMA 5.1 \it\ \ Let $\pi\colon \R^M\rightarrow \R^{M-1}$ 
        be the standard projection, 
        $\pi (x_1,\cdots ,x_M) = $  $(x_1,\cdots$ $ ,x_{M-1})$. Let $\Lambda$ be the lattice in
        $\R^{M-1}$ as in (4.1) and let $\eta >0$. 
        There is a constant $\omega >0$ such that for every $\tau >0$ the following holds. 
        Let $S\subset\R^{M-1}$ 
        be a simplex belonging to $\cD (\tau\Lambda)$. Suppose that $\psi$ is a Lipschitz 
        function in a neighbourhood of $S$ with Lipschitz 
        constant $\leq \omega $. Let $v_1, \cdots ,v_M$ be the vertices of $S$ and let  
        $w_1,\cdots ,w_M$ be the points in $\R^M$ given by 
        $w_j = (v_j,\psi (v_j)),\ 1\leq j\leq M$. Let $\Pi$ be the hyperplane in $\R^M$ 
        containing the points $w_1,\cdots ,w_M$ and let 
        $\Gamma $ be the sphere in $\Pi $, containing these points (that is, let $\Gamma $
        be the circumsphere of the $(M-1)$-simplex in $\Pi$ 
        with vertices $w_1,\cdots, w_M $). Then $\pi (\Gamma )$ is contained in the
        $(\tau\eta)$-neighbourhood of the circumsphere of the simplex $S$.
        \vskip 4mm
        \bf 6.\ Proof of Lemma 5.1 \rm
        \vskip 2mm Let $S\in \cD(\Lambda )$ and let $\eta >0$. 
        If we replace $\psi $ with $\psi + c $ where $c$ is a 
        constant, $\Pi$ will change to $\Pi + (0,c)$, \ 
        $\Gamma $ to $\Gamma + (0,c)$ and consequently $\pi (\Gamma ) $ will not change. Thus, 
        $\pi(\Gamma )$ remains unchanged if we subtract $\psi (v_M)$ from each $\psi (v_j),
        \ 1\leq j \leq M$. Thus,  
        $\pi (\Gamma )$ will be determined precisely 
        once we know $\beta_1 =\psi (v_1)-\psi (v_M), \cdots ,
        \beta_{M-1}= \psi (v_{M-1})-\psi (v_M)$. We shall show that $\pi  (\Gamma $) changes 
        continuously with $ (\beta_1,\ \cdots , \beta_{M-1})$ 
        near $(0,0,\cdots, 0)$ if $w_1=(v_1,\beta _1),\cdots , w_{M-1}=(v_{M-1}, \beta_{M-1})$ and 
        $ w_M= (v_M, 0)$. Note that when $\beta_1=\cdots =\beta _{M-1}= 0$, then 
        $\Gamma= \pi (\Gamma )$ is the circumsphere of
        $S$ in $\R^{M-1}$.  Let $w_0=(w_{01},\cdots w_{0, M-1}, 1)$ be the vector in $\R^M$ perpendicular to 
        $\Pi $ whose last component equals $1$. So $w_0$ must be perpendicular to $w_j-w_M,\ 1\leq j\leq M-1$, 
        so \ 
        $
        <w_j-w_M|w_0> = 0\ \ (1\leq j\leq M-1)
        $
        \ which, if $v_j=(v_{j1},\cdots ,v_{j, M-1}),\ \ 1\leq j \leq M-1$, is 
        the system of linear equations
        $$
        (v_{j1}-v_{M1})w_{01} +\cdots + (v_{j, M-1}-v_{M,M-1})w_{0,M-1} = - \beta_j \ \ (1\leq j\leq M-1).
        $$
        This is a system of $M-1$ linear equations for $M-1$ unknowns 
        $w_{01},\cdots ,w_{0, M-1} $ whose matrix is nonsingular, since, 
         $S$ being a $(M-1)$-simplex, the vectors 
         $v_j -v_M, \ 1\leq j\leq M-1$, are linearly 
         independent. Its solution 
         depends linearly on $(\beta_1,\cdots ,\beta _{M-1})$. When $\beta_1 =\cdots =\beta_{M-1} = 0$
         the solution 
         is the zero vector.
         In this case $w_0 = (0,\cdots ,0, 1)$. Let $z= (z_1,\cdots ,z_M)$ be 
         the center of the sphere in $\Pi$ that 
         contains $w_1,\cdots ,w_M$. Then $z$ is in $\Pi$ so
         $$
         <z-w_M|w_0> = 0 .
         \eqno (6.1)
         $$
         Further, for each $i, \ 1\leq i\leq M-1$,\ \ $z$ is at equal distance
         from $w_i$ and $w_M$, so $z$ is 
         contained in the hyperplane in $\R^M$ that passes through the
         midpoint of the segment joining $w_i$ and $w_M$, and
         is perpendicular to this segment, so $z$ must satisfy
         $$
         <[z-(w_i+w_M)/2] | [w_i-w_M]> = 0 .
         $$
         Thus,
         $$ 
         <z|[w_i-w_M]> =  (1/2) <[w_i+w_M]| [w_i-w_m]>\ \ (1\leq i\leq M-1).
         $$
         Together with (6.1) this becomes the following 
         system of linear equations for $z_1,\cdots, z_M$:
         $$ z_1(v_{i1}-v_{M1})+\cdots +z_{M-1}(v_{i, M-1}-v_{M, M-1})
         + z_M\beta_i =(|w_i|^2-|w_M|^2)/2 \ \ (1\leq i\leq M-1)
         $$
         $$
         z_1w_{01}+\cdots +z_{M-1}w_{0,M-1}+z_M = v_{M1}w_{01}+\cdots+ v_{M,M-1}w_{0, M-1}.
         $$
         Its matrix
         $$\left[ 
         \eqalign{
         &v_{11}-v_{M1}, \cdots , v_{1,M-1}-v_{M,M-1} , \beta_1 \cr
         &. . . . \cr
         &v_{M-1,1}-v_{M1}, \cdots , v_{M-1,M-1} , \beta_{M-1}\cr
         &w_{01}, \ \ \ \ \ \ \ \ \ \cdots \ \ \ \ \ \ \ \ \ \ \ \, w_{0, M-1} , \ 1\cr}\right]
         $$
         is nonsingular for $\beta_1=\cdots = \beta_{M-1}=0$ when $ w_{01}=\cdots = w_{0, M-1} = 0$.
         The matrix depends continuously on $(\beta_1, \cdots \beta_{M-1})$ and so do the
         right sides $(1/2)(|v_i|^2-|v_M|^2+\beta_i^2),\ 1\leq i,\leq M-1$, 
         and, since $w_0$ depends continuously on $(\beta_1,\cdots ,\beta_{M-1})$, 
         also $v_{M1}w_{01}+\cdots + v_{M,M-1}w_{0,M-1}$ depends continuously on 
         $(\beta_1,\cdots ,\beta_{M-1})$.
        \rm So the solution  $z= (z_1,\cdots, z_M) $, the center of the sphere $\Gamma$, depends continuously 
        on $(\beta_1,\cdots , \beta _{M-1})$ 
        near $(0,0,\cdots,0)$ and so does its radius $|z-w_M| = \bigl((z_1-v_1)^2+\cdots +(z_{M-1}-v_{M-1})^2+
        z_M^2)^{1/2}$. Recall that $\Pi$ passes through $w_M= (v_M, 0)$ and its perpendicular 
        direction $w_0$ changes continuously with 
        $(\beta_1,\cdots ,\beta_{M-1})$ so $\Pi$ changes continuously with
        $(\beta_1,\cdots ,\beta_{M-1})$. We have seen that the center 
        $z$ of the sphere $\Gamma $ in $\Pi$ and its radius also change continuously with 
        $(\beta_1,\cdots,\beta _{M-1})$ near $(0,0,\cdots 0)$. Thus, $\pi (\Gamma)$ changes 
        continuously with
        $(\beta_1,\cdots ,\beta_{M-1})$ near the origin where $\pi (\Gamma )=\Gamma$ is 
        the circumsphere of $S$ 
        when $\beta_1=\beta_2=
        \cdots\beta_{M-1} = 0$. Thus, $\pi (\Gamma )$ is contained
        in the $\eta$-neighbourhood of the circumsphere of
        of the simplex $S$ in $\R^{M-1}$ provided that 
        $\psi (v_1)-\psi (v_M),\cdots ,\psi (v_{M-1})-\psi (v_M)$ are small enough. 
        If $\psi $ is a Lipschitz function with the Lipschitz constant $\omega $ then 
        $|\psi (v_i)-\psi (v_M)|\leq \omega |v_i-v_M|,\ 1\leq i\leq M-1$, so 
        there is an $\omega $ such that if 
        $\psi $ is a Lipschitz function with the Lipschitz constant not 
        exceeding $\omega $ then $\pi (\Gamma )$ is contained 
        in the $\eta$-neighbourhood of 
        the circumsphere of the simplex $S$. Recall that every 
        simplex in 
        $\cD(\Lambda )$ is of the form $S_i+x, \ 1\leq i\leq \ell, 
        \ x\in \Lambda $. Repeating the reasoning above for each $S_i,\ 1\leq i\leq \ell$, we get the 
        Lipschitz constant that works for every simplex $S$ in $\cD(\Lambda )$. This completes the proof for $\tau =1$. 
        
        Now, let $\tau >0$ be arbitrary and let $S\subset \R^{M-1}$ be a simplex in $\cD (\tau \Lambda)$. 
        Let $\psi $ be a Lipschitz function with Lipschitz constant not exceeding 
        $\omega $ in a neighbourhood of 
        $S$, so its graph is given by $x_M= \psi (x_1,\cdots , x_{M-1})$. Introduce 
        new coordinates $X_1,\cdots ,X_M$ in $\R^M$ by 
        $x_j=    \tau X_j, \ 1\leq j\leq M$. In new coordinates we have $\tau X_M =
        \psi (\tau X_1,\cdots , \tau X_{M-1})$ so $X_M= \Psi (X_1,\cdots X_{M-1}) = 
        (1/\tau)\psi (\tau X_1,\cdots \tau X_{M-1})$. Both $\psi $ and $\Psi$ are 
        Lipschitz functions with the same Lipschitz constants, so in new coordinates 
        $\Psi $ is a Lipschitz function in a neighbourhood of $S$, which, in new 
        coordinates, belongs to $\cD(\Lambda )$. Thus, applying the first part of 
        the proof we see 
        that in new coordinates $\pi (\Gamma ) $ is
        contained in the $\eta$-neighbourhood of the circumsphere of $S$. In follows that in  
        old coordinates $\pi (\Gamma)$ is contained in the  $(\tau\eta)$-neighbourhood of the circumsphere 
        of $S$. This completes the proof.
        \vskip 4mm
        \bf 7.\ Polyhedral convex surface contained in a spherical shell\rm
        \vskip 2mm
        Let $\eta >0$ be as in Section 4 and let $\omega$ be the one given by Lemma 5.1.  Let again $W=U\times (1-\nu, 1+\nu)$ 
        where $\nu >0$ is small, $U$ 
        is a small open ball centered at the origin in $\R^{M-1}$ and let $R_0<1$ be
        so large that 
        $U\times\{ 1-\nu\}\subset R_0\B$. For every $r,\ R_0<r<1$,\ $W\cap b(r\B)
        =\{ (x,\psi_r(x)\colon\ x\in U\}$
        where the function $\psi_r $ is as in (5.1). We have $(\hbox{grad}(\psi_r))(x)=
        -(r^2-|x|^2)^{-1/2}x\ \ (x\in U)$ so 
        we may, passing to a smaller $U$ if necessary, 
        assume that $|(\hbox{grad} \psi_r)(x)|\leq\omega \ (x\in U,\ R_0<r<1)$ so that
        for each $r, \ R_0<r<1,\ \psi_r$ is a Lipschitz function on $U$
        with Lipschitz constant not exceeding $\omega $.
        
        Let $\Lambda$ be as in (4.1), let $\tau>0$ be small and let $R_0<r<1$. Let $\psi_r$ be as in (5.1). 
        Then $x\mapsto \Psi_r=(x,\psi_r (x))$ is a one to one map from $U$ onto $W\cap b(r\B)$.
        We now look at the points $\Psi_r(x),  x\in 
        (\tau\Lambda)\cap U$ and want to see them as vertices of a convex polyhedral 
        hypersurface in $\R^M$. 
        
        Consider a simplex $S\in\cD (\tau\Lambda)$ which is contained in $U$. Let $v_1,\cdots, v_M$ 
        be its vertices. We can extend the restriction of the function $\psi_r$ to this 
        set of vertices to a function $\varphi_r$ on all 
        $S$ by putting 
        $$
        \varphi _r \bigl(\sum_{j=1}^M \alpha_j v_j\bigr) = \sum_{j=1}^M \alpha_j\psi_r(v_j)\ \
        (0\leq \alpha_j\leq 1,\ 1\leq j\leq M,\ \sum_{j=1}^M \alpha_j = 1)
        $$
        to get an affine function $\varphi _r $ on $S$ so that 
        $x\mapsto \Phi _r(x) = (x,\varphi_r(x))$ is an affine map mapping $S$ to $\Phi_r(S)$, 
        the simplex with vertices 
        $\Psi_r (v_1),\cdots ,\Psi_r (v_M)$. We do this for every simplex $S\in\cD (\tau\Lambda)$ 
        that is contained in $U$. Thus, we get a piecewise linear function $\varphi_r$ on the union of 
        the simplices $S\in\cD (\tau\Lambda)$ contained in $U$ and so the union $C_r(\tau)$ of 
        all these $\Phi_r(S)$, 
        the graph of the function $\varphi_r$, is then a polyhedral surface in $\R^M$. We shall show 
        that the function $\varphi_r$ is 
        convex so that $C_r(\tau)$ is a convex polyhedral surface. Later we shall show that the part of 
        $C_r(\tau )$ contained in $W_0= U_0\cap (1-\nu, 1+\nu)$ with $U_0$ being a ball in $R^{M-1}$ 
        centered at the origin, strictly smaller than $U$, is a part of the boundary $bP$ of a suitable 
        convex polytope $P$.
        
        Given $S\in\cD (\tau\Lambda),\ S\subset U$, let $\Pi$ be the hyperplane in $\R^M$ that 
        contains $\Phi_r(S)$. Then $\Pi\cap b(r\B)$ is the sphere in $\Pi$ and which is the circumsphere
         of $\Phi_r(S)$ and which was denoted by $\Gamma $ in Section 5. By 
         Lemma 5.1,\ $\pi(\Gamma )$ is contained in the $(\tau\eta)$-neighbourhood of the circumsphere of 
         $S$ in $\R^{M-1}$. We know that the $\tau\eta$-neighbourhood of the closed ball 
         in $\R^{M-1}$ bounded by the circumsphere of $S$ contains no 
         other points of $\tau\Lambda$ than the vertices of $S$ which implies that all points of 
         $\Psi_r(U\cap(\tau\Lambda))$ other than the vertices of $\Phi_r(S)$ lie outside of the small
         "spherical cap" that $\Pi$ cuts out of $b(r\B)$, that is, outside of the 
         "small" part of $b(r\B)$ bounded by $\Gamma $. This shows that all other vertices of 
         the simplices in $C_r(\tau)$ that are not the vertices of $\Phi_r(S)$ are contained in the 
         \it open \rm halfspace of $\R^M$ bounded by $\Pi$ that contains the origin. Thus, 
         $\Phi_r(S)$ is a facet of $C_r(\tau )$. Since this holds for every $S\in\cD (\tau\Lambda), S\subset U$, it 
         follows that the surface $C_r(\tau )$ is convex.
         
         The simplices $\Phi_r(S)$ where $S\in\cD (\tau\Lambda),\ S\subset U$, have 
         all their vertices on $b(r\B)$. We want to estimate how far into $r\B$ they reach. 
         To do this, we need the following
         \vskip 2mm
         \noindent\bf PROPOSITION 7.1\ \it Let $0<r<1$, let $a\in b(r\B)$ and let $A\subset b(r\B)$ 
         be a set such that $|x-a|\leq \gamma$ for all $x\in A$ where $\gamma < r$.
         Then the convex hull of $A$ misses
         $\rho \overline\B$ where $\rho = r -\gamma^2/r$. \rm
         \vskip 1mm
         \noindent \bf Proof.\ \rm $A$ is contained in $\{ x\in b(r\B) \colon |x-a|\leq \gamma\} $. 
         With no loss of generality assume that $a = (r, 0, \cdots ,0)$. Then 
         $A\subset \{ x\in b(r\B) \colon \ (x_1-r)^2+x_2^2+\cdots +x_M^2\leq \gamma^2\}$ $\subset 
         \{ x\in b(r\B) \colon \ r^2-2x_1r+r^2\leq\gamma^2\} $ $= \{ x\in b(r\B) \colon \ 
         2r^2-2x_1r<2\gamma^2 \}= $ $\{ x\in b(r\B) \colon \ x_1> (r^2-\gamma^2)/r \}$ 
         $\subset \{ x\in r\overline \B \colon \ 
         x_1> (r^2-\gamma^2)/r\}$. The last set is a convex set that contains $A$ 
         and misses $\rho\overline\B$ which completes the proof.
         \vskip 1mm
         Denote by $d$ the length of the longest edge of simplices in 
         $\cD (\Lambda )$ so that $\tau d$ is the length of the longest edge of the simplices in 
         $\cD (\tau\Lambda)$. Since $\psi_r$ is a Lipschitz function with the Lipschitz 
         constant not exceeding $\omega $ the length of the longest edge 
         of the simplices $\Phi_r (S)$ where $S\in \cD (\tau\Lambda),\ S\subset U$, does not exceed
         $\sqrt{1+\omega^2}\tau d$. Now,  we use Proposition 7.1. If $R_0<r<1$ then $r-\gamma^2/r> 
         r-\gamma^2/R_0$. Thus, putting
         $$
         \lambda ={ {(1+\omega^2)d^2}\over {R_0}}
         $$
         we get the following 
         \vskip 2mm
         \noindent\bf PROPOSITION 7.2\ \it If $R_0<r<1$ then the simplices $\Phi_r (S)$ where 
         $S\subset \cD (\tau\Lambda),\ S\subset U, $ miss
         $\rho\overline\B$ where $\rho = r-\tau^2\lambda$. \rm
         \vskip 4mm
              \bf 8.\ A convex polytope with a prescribed part of the boundary\rm 
              \vskip 2mm
              We keep the meaning of $R_0, U, d $ and $\lambda $. Recall that 
              $U$ is an open ball in $\R^{M-1}$ centered at the origin. 
              Let $\mu $ be its radius. Let $0<\mu_0<\mu_1<\mu_2<\mu_3<\mu$ and let $U_i=\{ x\in\R^{M-1}\colon 
              |x|<\mu_i\}$,\ $W_i=U_i\times (1-\nu, 1+\nu),\ 0\leq i\leq 3$.
              
              Choose $\tau_0>0$ so small that 
              $$
              \tau_0 d < \hbox{min} \{ \mu -\mu_3, \mu_3 -\mu_2, \mu_2-\mu_1, \mu_1-\mu_0\} .
              \eqno (8.1)
              $$
              Then, since the maximal edge length of simplices in $\cD (\tau\Lambda)$ equals $\tau d$ it 
              follows that if $0<\tau<\tau_0 $ then
              
              - the simplices $S\in\cD (\tau\Lambda)$ that meet $U_0$ are contained in $U_1$
              
              - the simplices $S\in \cD (\tau\Lambda )$ that are contained in $U$  cover $U_2$. 
              \vskip 2mm
              \noindent\bf PROPOSITION 8.1\ \ \it There is a $\kappa >0$ such that whenever 
              $R_0\leq R\leq 1$ and 
              $R<R^\prime <R+\kappa$ then each hyperplane in $\R^M$ which meets $W_2
              \cap\bigl(R^\prime\overline\B
              \setminus R\overline\B\bigr)$ and misses $W_3\cap R\overline\B$ misses $R\overline \B$.
              
              \vskip 1mm
              \noindent\bf Proof \rm\ \ Suppose that there is no such $\kappa>0$. Then there are a 
              sequence $R_n ,\ R_0\leq R_n\leq 1\ (n\in\N)$, and a sequence $x_n\in W_2$, such that 
              $|x_n|>R_n\ (n\in\N)$ and such that $|x_n|-R_n\rightarrow 0$ as $n\rightarrow\infty$, 
              and for each $n$ a hyperplane $H_n$ through $x_n$ which misses $W_3\cap R_n\overline\B$
              and meets $R_n\overline\B\setminus W_3$. Since $|x_n|-\R_n\rightarrow 0$ as $n\rightarrow\infty$ 
              we may, passing to subsequences if necessary, with no loss of generality assume that 
              $R_n$ converges to an $R$ and $x_n$ converges to $x\in b(R\B)\cap \overline{W_2}$. 
              Since for each $n,\ H_n$\ misses $W_3\cap R_n\overline \B$ it follows that $H_n$ 
              converges to $H$, the hyperplane through $x$ tangent to $b(R\B)$ at $x$. 
              In particular, $H\cap (R\overline\B\setminus W_3)$ is empty, so for sufficiently large $n$, 
              $H_n\cap (R_n\overline\B\setminus W_3)$ must be empty, a contradiction.
              This completes the proof.
              
              \vskip 2mm
              With no loss of generality, passing to a smaller $\tau_0$ if necessary,
              we may assume that $\tau_0^2\lambda <\kappa$.
              Suppose now that $0<\tau<\tau_0$ and let $R_0\leq \rho<r<1$  
              where $\rho = r-\tau^2\lambda $.
              
              We know that the union $C_r(\tau)$ of the simplices $\Phi_r (S)$ 
              where $S\in\cD (\tau\Lambda), S\subset U, $ is a convex polyhedral 
              surface which, by Proposition 7.2, 
              is contained in $r\overline \B\setminus \rho\overline \B$. Each of 
              these simplices $\Phi_r(S)$ is contained in a hyperplane $H$. We want 
              that these hyperplanes miss 
              $\rho\overline\B$. Note that by (8.1) the simplices in $\cD(\tau\Lambda)$, contained 
              in $U$ cover $U_3$.  
             So the function $\varphi_r$ is well defined on $U_3$ and its graph 
             $C_r(\tau)\cap W_3$ 
             is contained in 
             $W_3\cap (r\overline\B\setminus\rho\overline\B)$. The function $\varphi_r$ 
             is piecewise linear and convex. Thus, 
             if $S\in\cD(\tau\Lambda)$ meets $U_2$ then, by (8.1),\ $S\subset U_3$ and by the 
             convexity of $\varphi_r$, the graph of $\varphi_r|U_3$ lies on one 
             side of the hyperplane 
             $H$ that contains $\Phi_r(S)$ which, in particular, implies that 
             $H$ misses $W_3\cap\rho\overline\B$ and thus, 
             by Proposition 8.1, $H$ misses $\rho\overline\B$. This shows that the part of
             $C_r(\tau )$ contained in 
             $W_2$  can be described in terms of the hyperplanes that miss $\rho\overline\B$. 
              So we find 
              $x_1,\cdots, x_n \in b\B$ and $\alpha_1, \cdots ,\alpha_n, \ \rho
              <\alpha_i\leq r\ (1\leq i\leq n)$, 
              such that 
              $$
              G_1=\{ x\in \overline\B\colon\ <x|x_i>\leq \alpha_i, 1\leq i\leq n\} .
              $$
              is a convex set containing $\rho\overline\B$ in its interior, and is such that 
              $W_2\cap bG_1 = W_2\cap C_r(\tau)$.
              \vskip 2mm
              \noindent\bf PROPOSITION 8.2\ \ \it Let $R_0<r<1$, let $0<\tau<\tau_0$, and 
              let $\rho=r-\tau^2\lambda > R_0$,  
              There is a convex polytope $P$ which contains 
              $\rho\overline \B$ in its interior, such that
              $bP\subset r\overline\B\setminus\rho\overline\B$, and such that every $\Phi_r (S)$ where 
              $S\in\cD (\tau\Lambda),\ S\subset U_1$, is a facet of $P$. \rm
              \vskip 2mm
              \noindent Proposition 8.2 implies in particular, that 
              $$
              W_0 \cap \hbox{skel}(P) = \Phi_r\bigl( U_0\cap \hbox{skel} (\cD (\tau\Lambda))\bigr)
              $$ 
              so that
              $$
              \pi \bigl(W_0\cap \hbox{skel}(P)\bigr) = U_0\cap \hbox{skel}\bigl(\cD (\tau\Lambda)\bigr) .
              $$
              \vskip 1mm
              \noindent\bf Proof.\ \rm To prove Proposition 8.2 we will find another convex set 
              $G_2$ whose 
              boundary outside $W_2$ will be a 
              polyhedral convex surface approximating $b(r\B)$ and such that
              $W_1\cap bG_2 = W_1\cap r\overline\B$ and then put $P= G_1\cap G_2$.
              To do this we first choose $\rho_1<r$ so close to $r$ that if $H$ is a hyperplane in $\R^M$ 
              passing through a point $x\in b(\rho_1\B)\setminus W_2$ tangent to $b(\rho_1\B)$ then 
              $H\cap W_1\cap r\overline B=\emptyset$. We will now use a finite number of these
              hyperplanes to 
              modify the part of $b(r\B)$ outside $W_1$ to get a convex polyhedral hypersurface
              contained 
              in $r\overline\B\setminus \rho_1\B$ which will be a part of $bG_2$. To do this, we need
              \vskip 2mm
              \noindent\bf PROPOSITION 8.3\ \it Let $x, y \in b\B$. Suppose that $ry$ is in the halfspace 
              $\{ z\in\R^M\colon\ <z|x>\leq \rho_1\}$, that is, in the halfspace bounded by the 
              hyperplane through $\rho_1x$, 
              tangent to $b(\rho_1\B)$ which contains the origin. Then $|x-y|\geq \sqrt{2(1-\rho_1/r)}$. \rm
              \vskip 1mm
              \noindent \bf Proof\ \rm Our assumption implies that $<ry|x>\leq \rho_1 $ so
              $<x|y>\leq \rho_1/r$ and so $|y-x|^2= 2-2<x|y>\geq 2-2\rho_1/r = 
              2(1-\rho_1/r)$ which completes the proof.
              \vskip 1mm
              \noindent Note that if $z\in b\B$ then $\{ y\colon\ <y|z>\leq \rho_1\}$ is 
              the halfspace bounded by the hyperplane through $\rho_1 z$ tangent to $b(\rho_1\B)$, 
              which contains the origin.
              \vskip 2mm
              \noindent\bf PROPOSITION 8.4\ \it Let $\cS$ be a subset of $b\B$. Let $0<\rho_1<r$ and let $0<
              \delta< \sqrt{ 2(1-\rho_1/r)}$. Assume that $z_1,\cdots , z_m \in \cS $ are such that
              $$
              \cS\subset \cup_{j=1}^m (z_j+\delta\B) .
              \eqno (8.2)
              $$
              Then the convex polyhedron 
              $$
              Q= \bigcap_{j=1}^m \{y\colon\ <y|z_j>\leq \rho_1\}
              $$
              does not meet $r\cS$. \rm
              \vskip 1mm
              \noindent \bf Proof\ \ \rm Suppose that $y\in\cS$ is 
              such that $ry\in Q$, that is $<ry|z_j> \leq \rho_1$ for 
              all $j,\ 1\leq j\leq m$. By Proposition 8.3 it follows
              that $|y-z_j|\geq \sqrt{2(1-\rho_1/r)} > \delta$
              for all $j, 1\leq j\leq m$, 
              which contradicts (8.2). This completes the proof. 
              \vskip 2mm
              \noindent We now proceed to finish the proof of Proposition 8.2. 
              Let $\cT = b(r\B)\setminus W_2$. Choose 
              $\delta, \ 0<\delta<\sqrt{2(1-\rho_1/r)}$, and then choose $z_1, \cdots ,z_m\in b\B$ such that
              $$
              {1\over r}\cT \subset \cup_{j=1}^m (z_j+\delta\B) .
              $$
              Set
              $$
              G_2=\{ y\in r\overline\B\colon\ <y|z_j>\leq \rho_1 \ (1\leq j\leq m)\}
              $$
              and let \ $P=G_1\cap G_2$, so 
              $$
              P= \{ x\in\overline\B\colon <x|x_i>\leq \alpha_i, 1\leq i\leq n, \ 
              <x|z_j> \leq \rho_1, \ 1\leq j\leq m\} .
              $$
              By construction, $P$ contains $\rho\overline\B$ in its interior. Moreover, it is easy to see that 
              $$
               P= \{ x\in\R^M\colon \ <x|x_i>\leq \alpha_i, 1\leq i\leq n, \ 
               <x|z_j> \leq \rho_1, \ 1\leq j\leq m\} ,
              $$
              so $P$ is a convex polytope contained in $r\overline\B$ and, by construction, is such that every 
              $\Phi_r (S)$ where 
              $S\in\cD (\tau\Lambda),\ S\subset U_1$, is a facet of $P$. Proposition 8.2 is proved. 
              
              It is clear that all we have done so far will work in the same way for any lattice 
              $\tau (\Lambda +q)$. Summing up what we have proved so 
              far we get our main Lemma 8.5. Recall that
              $\pi (z_1,\cdots ,z_M)= (z_1,\cdots ,z_{M-1})$.
              \vskip 2mm
              \noindent\bf LEMMA 8.5\ \it There are $R_0,\  0<R<1,\  \nu >0, \ \tau_0 >0, \ 
              \lambda >0$ 
              and a small open ball $U_0 \subset R^{M-1}$ centered
              at the origin, such that $U_0\times\{ 1-\nu\} 
              \subset R_0\B$ and such that if $W_0=U_0\times (1-\nu, 1+\nu)$ then the following holds:\ 
              
              \noindent For each $\tau, \ 0<\tau <\tau_0$, for each $r$ such that 
              $$
              R_0<r-\lambda\tau^2<r<1
              $$
              and for each $q\in\R^{M-1}$ there is a convex polytope $P$ 
              contained in $r\overline\B$ and containing 
              $(r-\lambda\tau^2)\overline\B$ in its interior and such that
              $\pi $ maps $W_0\cap\hbox{skel}(P)$ onto $U_0\cap\hbox{skel}(\cD(\tau(\Lambda + q)))$. \rm
              \vskip 4mm
         \bf 9. Small blocks of convex polytopes \rm
         \vskip 2mm
         Let $\Lambda $ be as in (4.1) and let $E(\Lambda )$ be the fundamental parallelotope 
         for $\Lambda$, that is,
         $$
         E(\Lambda ) = \{\theta_1e_1+\cdots +\theta_{M-1}e_{M-1}\colon\ 0\leq \theta_i <1,\ 
         1\leq i\leq M-1\}.
         $$ 
         Given $q\in\R^{M-1}$ define $\cS (q)= \hbox{skel} (\cD(\Lambda + q))$. 
         Clearly $\cS (q)=\cS (0)+q$. 
         Recall that all our tessellations are periodic so 
         $$
         \cS (q) +\sum_{j=1}^{M-1}n_je_j = \cS (q)
         $$ 
         for every $q\in \R^{M-1}$ and every $n_j\in Z,\ 1\leq j\leq M-1$. Thus, if $w\in \cS (q_1)\cap 
         \cS (q_2)$ there are $n_j,\ 1\leq j\leq M-1$ such that if $w_0= w- \sum_{j=1}^{M-1}n_je_j 
         \in E(\Lambda) $ then $w_0\in E(\Lambda) \cap\cS (q_1)\cap \cS (q_2)$. Thus, if 
         $\cS (0)\cap\cS (q_1)\cap\cdots \cap \cS(q_{M-1})\cap E(\Lambda)=\emptyset$ then 
         $\cS (0)\cap\cS (q_1)\cap\cdots \cap \cS(q_{M-1}) =\emptyset$
         \vskip 2mm
         \noindent \bf PROPOSITION 9.1\ \it Given $\varepsilon >0$ there 
         are $q_1,\cdots, q_{M-1}$, 
         $|q_i|<\varepsilon,\ 1\leq i\leq M-1$, 
         such that $\cS (0)\cap \cS(q_1)\cap\cdots\cap \cS (q_{M-1})= \emptyset$ . \rm 
         \vskip 2mm
         \noindent We need the following 
         \vskip 2mm
         \noindent\bf PROPOSITION 9.2\ \it Let $H$ be a hyperplane in $\R^{M-1}$. Let $\tilde H$ 
         be the hyperplane in $\R^{M-1}$ parallel to $H$ which passes through the origin and assume 
         that $q\in\R^{M-1}$,\ $q\not\in\tilde H$. Let $L$ be a $k$-plane in $\R^{M-1}$ where 
         $1\leq k\leq M-2$. Then 
         either $L\subset H+tq $ for some $t\in\R$ or else $L$ intersects $H+tq$ transversely 
         for every $t\in\R$. \rm
         \vskip 2mm
         \noindent\bf Proof\ \ \rm Obvious.
         \vskip 2mm
         \noindent We shall say that a $k$-plane $L$ is transverse 
         to a hyperplane $G$ if it is not contained in $G$. In this case either $L$ misses 
         $G$ or else $L$ intersects 
         $G$ transversely (and $L\cap G$ is a $(k-1)$-plane). So the proposition says that $L$ is 
         transverse to the hyperplane $H+tq$ for each $t$ except for perhaps one value of $t$.
         \vskip 1mm
         \noindent\bf Proof of Proposition 9.1\ \rm Take a large ball $B$ centered at the 
         origin and consider the family of all those hyperplanes that contain a facet 
         of a simplex $S\in\cD(\Lambda )$ contained in $B$. There are finitely many of these 
         hyperplanes, denote them by $L_1,\cdots, L_p$ and their union by $\L$. For 
         each $j,\ 1\leq j\leq p$, let $\tilde L_j$ be the hyperplane parallel to $L_j$ passing through the 
         origin. Choose $q\in\R^{M-1}$ so that $q$ belongs to no $\tilde L_j,\ 1\leq j\leq p$. 
         Let $
         \varepsilon >0$. By the dicussion at the beginning of this section the proposition 
         will be proved once we have proved that there are $t_j,\ \varepsilon>t_1>\cdots >t_{M-1}>0$ 
         such that
         $$
         \L\cap (\L+t_1q)\cap\cdots\cap (\L+t_{M-1}q) = \emptyset 
         $$
         and then we put $q_j=t_jq,\ 1\leq j\leq M-1$.
         
         By Proposition 9.2 for each $j,\ 1\leq j\leq p$, and for
         each $t,\ 0<t<\varepsilon $, except perhaps finitely many, $L_j+tq$ 
         is transverse to each $L_k,\ 1\leq k\leq p$. So there is a $t_1, \ 0<t_1<\varepsilon $, that 
         works for all $L_j,\ 1\leq j\leq p$, so that $\L\cap (\L +t_1q)$ is a union of finitely 
         many $(M-3)$-planes. Suppose that $1\leq\ell\leq M-3$ and suppose that we have found $t_1,\cdots, 
         t_\ell$, 
         $\varepsilon>t_1> t_2 >, \cdots >t_\ell>0$, such that $\L\cap(\L+t_1q)\cap\cdots\cap 
         (\L+t_\ell q)$ is 
         a finite union of $(M-2-\ell )$-planes . Applying Proposition 9.2 we find
         $t_{\ell +1},\ 0<t_{\ell +1}<t_\ell$, such that 
         $\L\cap (\L+t_1q)\cap\cdots\cap (\L+t_{\ell+1}q)$ is a 
         finite union of $(M-3-\ell)$-planes. Thus, step by step we arrive to the point where $\L\cap (\L+t_1q)\cap
         \cdots \cap (\L+t_{M-2}q)$ is a finite set of points whose intersection with $\L+t_{M-1}q$ with a suitable 
         chosen $t_{M-1}, \ 0<t_{M-1}<t_{M-2}$ is empty. This completes the proof. 
         \vskip 2mm
         \noindent\bf LEMMA 9.3\ \it Let $q_0=0$ and let $q_1,\cdots ,q_{M-1}$ be as in Proposition 9.1. Let
         $$
         \cS_i = \hbox{\rm skel\it}\bigl(\cD (\Lambda +q_i)\bigr)\ \ (0\leq i\leq M-1) .
         $$ 
         There is a $\mu >0$ such that whenever $x_i\in \cS_i,\ 0\leq i\leq M-1$, we have 
         $$
         |x_1-x_0|+|x_2-x_1| +\cdots + |x_{M-1} -x_{M-2}| \geq \mu .
         \eqno (9.1)
         $$
         \rm
         \noindent\bf Proof\ \ \rm Assume that there is 
         no $\mu >0$ such that (9.1) holds  whenever $x_i\in \cS_i,\ 0\leq i\leq M-1$. 
         Then there are sequences 
         $x_{i,n}\in\cS_i,\ 0\leq i\leq M-1,\ n\in\N$ such that 
         $$
          |x_{1n}-x_{0,n}|+|x_{2n}-x_{1n}| +\cdots + |x_{M-1,n} -x_{M-2,n}| 
          \eqno (9.2)
         $$
         tends to zero as $n\rightarrow\infty$.
         Since $\cS_i $ are periodic, that is, 
         $$
         \cS_i = \cS_i + \sum_{k=1}^{M-1} m_ke_k\ \ (0\leq i\leq M-1)
         $$
         whenever $m_k\in Z,\ 1\leq k\leq M-1$, adding for each $n$, 
         a suitable $\sum_{k=1}^{M-1}m_{k,n}e_k$ 
         to all $x_{0n}, x_{1n},\cdots , x_{M-1,n}$ 
         where $m_{k,n}\in Z, 1\leq k\leq M-1$ - note that doing this,
         the sum (9.2) remains unchanged - we may,
         with no loss of generality assume that $x_{0n}\in E(\Lambda)$ for all $n$, so, by compactness, 
         we may, after 
         passing to a subsequence if necessary, assume that $x_{0n}$ converges to some $x_0$.
         Since $\cS_0$ is closed,\ $x_0\in \cS_0$. Since (9.2) tends to zero as $n\rightarrow\infty$ it follows 
         that for each $j,\ 0\leq j\leq M-1$, the sequence $x_{jn}\in\cS_j$ converges to the
         same limit $x_0$ which must be in $\cS_j$ since $\cS_j$ is closed. Thus, $x_0$ is 
         contained in the intersection $\cS_0\cap\cdots\cap\cS_{M-1}$ contradicting the fact
         that this intersection is empty. This completes the proof. 
         \vskip 1mm
         Let $q_i,\ 0\leq i\leq M-1$ be as in Lemma 9.3. For each $\tau>0$ we have 
         $$
         \hbox{skel}\bigl(\cD (\tau(\Lambda + q))\bigr) = \tau \hbox {skel}\bigl(\cD(\Lambda + q)\bigr)
         $$
         so by Lemma 9.3 it follows that if $\tau >0$,  and if 
         $x_i\in\hbox{skel}\bigl(\cD 
         (\tau(\Lambda + q_i))\bigr)$, $0\leq i\leq M-1$ then
         $$
         |x_1-x_0|+|x_2-x_1| +\cdots + |x_{M-1} -x_{M-2}| \geq \tau\mu .
         $$ 
         \vskip 2mm
         \noindent\bf LEMMA 9.4\ \it Let $0<\tau<\tau_0$ and suppose that 
         $$
         R_0<r-M\tau^2\lambda <r<1 .
         $$
         There are convex polytopes $Q_j,\ 0\leq j\leq M-1$, such that
         $$
         ((r-M\tau^2\lambda)\overline\B \subset \hbox{\rm Int}Q_0\subset\hbox{\rm Int}Q_1\subset\cdots\subset Q_{M-1}\subset r\overline\B
         $$
         such that for each $j,\ 0\leq j\leq M-1$, 
         $$
         \pi (W_0\cap\hbox{\rm skel}(Q_j)) = U_0\cap\hbox{\rm skel}(\cD(\tau(\Lambda + q_j)).
         $$
         Thus, 
         $$
         \left.
         \eqalign{
         &\hbox{if\ \ } x_j\in W_0\cap\hbox{skel}(Q_j)\ \ (0\leq j\leq M-1)\hbox{\ \ then}\cr
         & |x_1-x_0|+\cdots +|x_{M-1}-x_{M-2}| \geq \tau\mu\cr}
         \right\} 
         \eqno (9.3)
         $$
         \noindent \bf Proof\ \rm Let $0\leq j\leq M-1$.  By Lemma 8.5 there is a convex polytope 
         $Q_j$ containing 
         $\bigl( r -(M-j)\tau^2\lambda \bigr)\overline\B$ in its interior and contained in 
         $\bigl( r -(M-(j+1))\tau^2\lambda \bigr)\overline\B$ such that $\pi $ maps 
         $W_0\cap\hbox{skel}(Q_j)$ onto 
         $U_0\cap \hbox{skel}\bigl(\cD(\tau(\Lambda + q_j))\bigr).$ Thus, if $x_j,\ 0\leq j\leq M-1$, are 
         as in (9.3) then $\pi (x_j)\in\hbox{skel}\bigl(\cD(\tau(\Lambda+q_j))\bigr)\ \ (0\leq j\leq M-1)$ 
         and hence by the discussion preceding Lemma 9.4 we have
         $$
         |\pi(x_1)-\pi (x_0)| +\cdots +|\pi(x_{M-1}-\pi (x_{M-2})|\geq \tau\mu
         $$
         so
         $$
         |x_1-x_0|+\cdots+|x_{M-1}-x_{M-2}| \geq \tau\mu .
         $$
         This completes the proof. 
         \vskip 1mm
         We shall call the family $\{ Q_0, Q_1, \cdots Q_{M-1}\}$ as above a \it small 
         block of convex polytopes \rm with boundaries contained in 
         $r\overline \B\setminus (r- M\tau^2\lambda)\overline \B$. More generally if $A\colon\ 
         \R^M\rightarrow\R^M$ is a rotation, that is, $A\in SO(M)$,  then we will call the family 
         $\{ A(Q_0), A(Q_1), \cdots ,
         A(Q_{M-1})\}$ also a small block of convex polytopes. 
         \vskip 4mm
         \bf 10.\ Large blocks of convex polytopes\ \rm
         \vskip 2mm
         In previous section we constructed a small block of convex polytopes, that is, given 
         $\rho, \ R_0<\rho-M\tau^2\lambda <\rho<1$, we constructed
         convex polytopes $Q_j,\ 0\leq j\leq M-1$, such that
         $$
         (\rho-M\tau^2\lambda)\overline\B \subset \hbox{Int}Q_0\subset Q_0\subset\cdots\subset 
         \hbox{Int}Q_{M-1}\subset Q_{M-1}\subset\rho\overline\B ,
         $$
         and such that (9.3) holds. An analogous statement holds if 
         we apply a rotation $A$ to all polytopes $Q_j,\ 1\leq j\leq M-1$, to get a new small block of 
         convex polytopes $R_j= A(Q_j),\ 0\leq j\leq M-1$, which have the property that 
         if $x_j\in A(W_0)\cap\hbox{skel}(R_j)\ (0\leq j\leq M-1)$ then 
         $$
         |x_1-x_0|+\cdots +|x_{M-1}-x_{M-2}| \geq \tau\mu .
         $$
         It is perhaps appropriate to mention that different 
         convex polytopes $Q^\prime$ and $Q^{\prime\prime}$ in the family of convex 
         polytopes that we are constructing have always their boundaries in disjoint 
         spherical shells so that if $Q^\prime\subset \hbox{Int}Q^{\prime\prime}$ and if $A$ 
         is a rotation then $A(Q^\prime)\subset \hbox{Int}Q^{\prime\prime}$.
         
         We now choose rotations $A_1=\hbox{Id}, A_2,\cdots ,A_L$ so that the open sets
         $$
        W_{0j}=A_j(W_0),\ 1\leq j\leq L,\ \hbox{cover\ } b\B,\hbox{\ that is,\ }
        b\B\subset \bigcup_{j=1}^L W_{0j} .
        \eqno (10.1)
        $$
        We now construct what we call a \it large block of 
        convex polytopes \rm that will have a property analogous to (9.3) for a sequence 
        $x_j,\ 0\leq j\leq M-1$  
        contained in any of the sets $W_{0j},\ 1\leq j\leq L$. Roughly speaking, we shall take 
        $\rho_0<\rho_1<\cdots <\rho_L$ and for each spherical shell 
        $\cS_k = \rho_k\overline\B\setminus\rho_{k-1}\overline\B,\ \ 1\leq k\leq L$, we 
        shall construct a small block $\cB_k$ of convex polytopes with boundaries contained in $\cS_k$ which has 
        the property (9.3) for $Q_j\in\cB_k,\ 0\leq j\leq M-1$. 
        Then we will rotate each $\cB_k$ by
        $A_k$, to form an $L$-tuple of smal blocks $A_1(\cB_1), A_2(\cB_2),\cdots, A_L(\cB_L)$, 
        and then arrange all the convex polytopes of these $A_j(\cB_j) $ into a single sequence, 
        Here is the exact formulation.
        \vskip 2mm
        \noindent\bf LEMMA 10.1\ \it Given $\tau,\ 0<\tau<\tau_0$, and $r$ such that
        $$
        R_0<r-ML\tau^2\lambda<r<1
        $$
        there is a family of convex polytopes $C_j,\ 0\leq j\leq ML-1$, such that
        $$ 
        (r-ML\tau^2\lambda )\overline\B    \subset \hbox{\rm Int}C_0\subset C_0\subset
        \hbox{\rm Int}C_1\subset
        \cdots\subset \hbox{\rm Int}C_{ML-1}\subset C_{ML-1}\subset r\overline\B
        $$
        which has the property that if \ \ $1\leq k\leq L$, and if $x_j\in W_{0k}\cap\hbox{\rm skel}C_j,\ 
        0\leq j\leq ML-1$, \ then 
        $$
        |x_1-x_0|+|x_2-x_1|+\cdots + |x_{ML-1}-x_{ML-2}| \geq \tau\mu .
        $$ \rm
        \vskip 2mm
        \noindent We shall call the family $\cC =\{C_0,C_1,\cdots , C_{ML-1}\} $ as above 
        \it a large block of convex 
        polytopes \rm with boundaries contained in
        $r\overline\B\setminus (r-ML\tau^2\lambda)\overline\B$.
        \vskip 2mm
        \noindent\bf Proof\ \rm Let
        $$
        \rho_j = r-M(L-j)\tau^2\lambda\ \ (0\leq j\leq L).
        $$
        For each $j,\ 1\leq j\leq L$, there is a small block $\cB_j$ of convex polytopes 
        with boundaries contained in 
        $\rho_j\overline\B\setminus\rho_{j-1}\overline\B$ such that (9.3) holds.
        
        Let $A_j,\ 1\leq j\leq L$, be rotations of $\R^M$ satisfying (10.1). 
        For each $j,\ 1\leq j\leq L $, 
        form a new small block $\cA_j = \{ A_j(P)\colon\ P\in\cB_j\} = \{
        C_{j0}, C_{j1},\cdots, C_{j, M-1}\} $ where
        $$
        \rho_{j-1}\overline\B\subset\hbox{Int}(C_{j0})\subset\hbox{Int}(C_{j1})\subset
        \cdots\subset C_{j,M-1}\subset \rho_j\overline \B
        $$
        such that if 
        $$
        x_i\in W_{0k}\cap\hbox{skel}(C_{ki})\ , \ \ 0\leq i\leq M-1,
        $$
        then 
        $$
        |x_1-x_0|+\cdots +|x_{M-1}-x_{M-2}| \geq \tau\mu .
        $$ 
        Now, write all $C_{ji}$ into a single sequence $C_{10}, C_{11},\cdots C_{1, M-1},
        C_{20},\cdots, 
        C_{2, M-1},\cdots ,$ $C_{L0}$, $C_{L1},\cdots, $ $C_{L, M-1}$, in other words
        $$ 
        C_{(j-1)M+i} = C_{ji}\ \ (1\leq j\leq L,\ 0\leq i\leq M-1).
        $$
        It is easy to see that the convex polytopes $C_0,\ C_1, \cdots\ ,C_{LM-1}$ have 
        all the required properties. This completes the proof. 
        \vskip 4mm
        \bf 11.\ Completion of the proof of Theorem 2.1 and the proof of Corollary 2.2\ \rm
        \vskip 2mm
        We keep the meaning of $R_0$ and $\tau_0$. Recall that by (10.1) the open sets $W_{0 j}= A_j(W_0),\ 1\leq j\leq L$, cover $b\B$. Thus
        $$
        \left.\eqalign
        {
        &\hbox{if\ \ }x_n\in\B\ \ \hbox{converges to \  } x\in b\B \hbox{\ then there are \ \ }n_0 \cr
        &\hbox{and\ \ } j,\ 1\leq j\leq L, \hbox{\ such that\ \ }x_n\in W_{0j}\ \ (n\geq n_0)\cr} 
        \right\} \eqno (11.1)
        $$
        To complete the proof of Theorem 2.1 we shall construct a sequence $r_j,\ 
        R_0<r_1<\cdots <r_j<\cdots <1$, 
        converging to $1$, and for each $j\in\N$ we shall construct a large block $\cC_j = \{ C_{j0}, C_{j1}, \cdots ,
        C_{j, LM-1}\} $ of convex polytopes such that
        $$
        r_j\overline\B\subset\hbox{Int}C_{j0}\subset C_{j0}\subset \hbox{Int }C_{j1}\subset
        \cdots\subset \hbox{Int}C_{j, LM-1}\subset C_{j,LM-1}\subset r_{j+1}\overline\B
        \eqno (11.2)
        $$
        such that writing all polytopes of all large blocks into a single sequence, that is,
        $$
        P_{(j-1)LM+k}= C_{jk}\ \ (0\leq k\leq LM-1, j\in\N )
        \eqno (11.3)
        $$
        we get our sequence $P_n$ of convex polytopes with the desired properties.
        
        \noindent To do this, choose $r_1,\ R_0<r_1<1$, and a decreasing sequence of positive numbers 
        $\tau_j, \ \tau_1<\tau_0$, such that
        $$
        \sum_{j=1}^\infty \tau_j^2 = {{1-r_1}\over{ML\lambda}} \hbox{\ \ and such that\ \ } 
        \sum_{j=1}^\infty \tau_j\hbox{\ \  diverges},
        \eqno (11.4)
        $$
        and then let
        $r_{j+1}=r_j+ML\tau_j^2\lambda\ \ (j\in\N).
        $
        Note that the equality in (11.4) means that the sequence $r_j$ converges to $1$ as
        $j\rightarrow\infty$. 
        
        Use Lemma 10.1 to show that for each $j\in\N$ there is a large block 
        $$
        \cC_j =\{ C_{j0}, C_{j1},\cdots C_{j, LM-1}\} 
        $$
        of convex polytopes satisfying (11.2) and having the property that 
        $$\left.\eqalign
        {&\hbox{if for some\ }k,\ 
        1\leq k\leq L, \ \ x_\ell\in W_{0k}\cap\hbox{\rm skel}(C_{j\ell})
        \hbox{\ \ for each\ } \ell, 
        \ 0\leq\ell\leq LM-1\ \cr 
        &\hbox{then\ \ } |x_1-x_0|+|x_2-x_1|+\cdots+|x_{LM-1}-x_{LM-2}|\geq 
        \tau_j\mu
        \cr}\right\}
        \eqno (11.5)
        $$
        Define the sequence $P_n$ of convex polytopes by 
        writing all polytopes $C_{jk}$ into a single sequence as in (11.3). Obviously
        $$
        P_0\subset \hbox{Int}P_1\subset P_1\subset\cdots\subset \B,\ \ \bigcup_{j=0}^\infty P_j= \B .
        $$
        Now, let $w_n\in\hbox{skel}(P_n)\ (n\in\N)$. To complete the proof of Theorem 2.1 we 
        must show (2.2). We know that it is enough to show this for sequences $w_n$ that converge. 
        So assume that $w_n$ converges. The properties of $P_n$ imply that the limit of the 
        sequence $w_n$ is contained in $b\B$. By (11.1) 
        there are $k, 1\leq k\leq L$, and $n_0$ such that $w_n\in W_{0k}\ \ (n\geq n_0)$. Let 
        $j_0$ be so large that $j_0 NL\geq n_0$. By (11.5), for each $j\geq j_0$, the large 
        block of polytopes $\cC_j$ adds at least $\tau_j\mu$ to the sum of the absolute values of 
        differences of consequtive $w_j$-s, that is,  for each $j\geq j_0$ we have
        $$
        |w_{(j-1)ML+1}-w_{(j-1)ML}| +\cdots +|w_{jML-1}-w_{jML-2}| \geq \tau_j\mu .
        $$
        It follows that for each $j\geq j_0$ there is a $N(j)<\infty$ such that 
        $$
        \sum _{i=1}^{N(j)}|w_i-w_{i-1}| \geq \sum_{k=j_0}^j\tau_k\mu .
        $$ 
        The fact that the series $\sum_{i=1}^\infty \tau_j $ diverges implies (2.2). 
        The proof of Theorem 2.1 is complete.

        \vskip 2mm
        \noindent\bf Proof of Corollary 2.2 \rm  Let $p\colon \ [0,1)\rightarrow \B$ 
        be a path such that $|p(t)|\rightarrow 1$ as $t\rightarrow 1$, 
        and such that for all sufficiently large $n\in\N$,\ \ 
        $p([0,1))$ meets $bP_n$ only at $\cU_n$.   
        Since $|p(t)|\rightarrow 1$ as $t\rightarrow 1$, it follows that $p(t)$ has to leave each
        $P_n$ so there are an $n_0$ and a sequence $t_j$,
        $$
        t_{n_0}<t_{n_0+1}<\cdots <1,\ \ \lim_{n\rightarrow\infty}t_n=1,
        $$
        such that $p(t_n)\in bP_n$ for each $n\geq n_0$. Thus, by our assumption, 
        passing to a larger $n_0$ if necessary, 
        we may assume that  $p(t_n)\in\cU_n$ for each $n\geq n_0$. 
        Thus, for each $n\geq n_0$ there is an $x_n\in \hbox{skel}(P_n)$  such 
        that $|x_n-p(t_n)|<\theta_n$. 
         For $n\geq n_0$ we have $|p(t_{n+1})-p(t_n)| 
        \geq|x_{n+1}-x_n|-|p(t_{n+1}-x_{n+1}| - |p(t_n)-x_n| \geq |x_{n+1}-x_n|-\theta_{n+1} - \theta_n$.
        It follows that
        $$
        \sum_{n=n_0}^\infty |p(t_{n+1})-p(t_n)|\ \geq  \sum_{n=n_0}^\infty |x_{n+1}-x_n| - 2
        \sum_{n=n_0}^\infty \theta_n.
        $$
        Since, by Theorem 2.1,  the series $\sum_{n=n_0}^\infty |x_{n+1}-x_n| $ diverges and since the series 
        $\sum_{n=n_0}^\infty \theta_n$  converges it follows that the series 
        $$
        \sum_{n=n_0}^\infty |p(t_{n+1})-p(t_n)|
        \eqno (11.6)
        $$ 
        diverges. Since the sequence $t_m$ 
        increases it follows that the length of $p([t_{n_0},1))$ is bounded from below  by 
        the sum of the series (11.6). Since this series diverges it follows that $p$ 
        has infinite length. 
        This completes the proof of Corollary 2.2.   
        \vskip 4mm
        \bf 12.\ Proof of Theorem 1.1 \rm
        \vskip 2mm
        As we know, every convex polytope $P\subset \R^M$ which contains the origin 
        in its interior can be written as 
        $$
        P = \bigcap _{i=1}^n K(x_i,1) = \bigcap _{i=1}^n\{ y\in\R^M\colon\ <y|x_i>\leq 1\} 
        \eqno (12.1)
        $$
        with $x_i\in\R^M\setminus\{ 0\}\  ,\ 1\leq i\leq n$. We assume 
        that the representation (12.1) is irreducible, so
        $$
        bP=\bigcup_{i=1}^nH(x_i,1)\cap P = \bigcup_{i=1}^n \{y\in \R^M\colon\ <y|x_i>=1\}\cap P ,
        $$
        and the sets $F_j=H(x_j,1)\cap P,\ 1\leq j\leq n$, are precisely the facets of $P$. Recall that 
        $\hbox{skel}(P) = \bigcup_{i=1}^n[F_i\setminus\hbox{ri}(F_i)]$.
        \vskip 2mm
        \noindent\bf PROPOSITION 12.1\ \ \it Let $P$ be as above. Let $\theta >0$.
        There is an $\eta >0$ such that for each $i,\ 1\leq i\leq n$, the set
        $$
        bP\cap \{ y\in\R^M\colon\ 1-\eta < <y|x_i><1\} 
        $$
        is contained in the $\theta$-neighbourhood of $\hbox{\rm skel}(P)$ in $bP$.
        \vskip 2mm 
        \noindent\bf Proof\ \ \rm Assume that Proposition 12.1 does not hold so that there 
        are $i,\ 1\leq i\leq n$, and 
        $\theta >0$ such that for each $\eta >0$ there is some $y\in bP$ such that 
        $1-\eta < <y|x_i> < 1$ and 
        $\hbox{dist}(y, \hbox{skel}(P))\geq\theta$. So there is a sequence $y_n\in bP$ 
        such that $<y_n|x_i> < 1 \ (n\in \N),\ $\ $<y_n|x_i>
        \rightarrow 1$\ as\ $n\rightarrow \infty$ and such that $\hbox{dist}(y_n , 
        \hbox{skel}(P)\geq \theta $ for all $n$. By compactness we may, after passing to a 
        subsequence if necessary, assume that $y_n$ 
        converges to  $y_0\in bP$. Clearly $y_0\in H(x_i,1)$. Since $y_0\in bP$ it follows 
        that $y_0$ belongs to the facet 
        $F_i = P\cap H(x_i, 1)$. Since $\hbox{dist}(y_0, \hbox{skel} (P)) 
        \geq \theta$ it follows that $y_0\in\hbox{ri}(F_i)$. 
        On the other hand, since $y_n\in bP\setminus F_i$ it follows that $y_n\in
        \cup_{j=1, j\not=i} F_j$. Passing to a subsequence if necessary we may assume 
        that there is a $j\not=i$, such 
        that $y_n\in F_j$ for all $n$. Since $F_j$ is closed it follows that $y_0\in F_j$. 
        Thus $y_0$, a relative interior point of
        the facet $F_i$, belongs to a different facet $F_j$ which is impossible. This completes the proof.
        \vskip 2mm
        \noindent\bf REMARK\ \ \rm Note that if $\cU$ is the $\theta$-neighbourhood of 
        $\hbox{skel}(P)$ and if $\eta$ is as above then for each $j,\ 1\leq j\leq n$, the set 
        $\{ y\in \R^M\colon\ <y|x_j>\leq 1-\eta\}$ contains $\cup_{i=1, i\not= j}^n[F_i\setminus \cU]$.
        \vskip 2mm
        We now move to $\C^N = \R^{2N}$ and denote by $<|>$ the Hermitian inner product in $\C^N$. Note that 
        $\Re(<|>) $ is then the standard inner product in $\R^{2N}$.
        \vskip 2mm
        \noindent\bf LEMMA 12.2\ \ \it Let $P$ be a convex polytope in $\C^N$
        and let $K\subset \hbox{\rm Int}(P)$ be 
        a compact set. Let $\theta >0$ and let \ $\cU\subset bP$ be 
        the $\theta$-neighbourhood of $\hbox{skel}(P)$ 
        in $bP$.
        Given $\varepsilon >0$ and $L<\infty$ there is a polynomial 
        $f\colon\ \C^N\rightarrow \C$ such that 
        $$
        \ \ \Re (f(z)) \geq L\ \ (z\in bP\setminus \cU) \hbox{\ \ and\ \ } |f(z)|<\varepsilon\ \ (z\in K).  
        $$\rm
        \vskip 1mm
        \noindent\bf Proof\ \ \rm With no loss of generality assume that the origin is an interior point of $P$. 
        There are $n\in\N$ and
        $w_1, w_2, \cdots , w_n \in \C^N\setminus \{ 0\}$ such that 
        $$
        P=\bigcap_{i=1}^n \{ z\in\C^N\colon\ \Re (<z|w_i>) \leq 1 \} 
        \eqno (12.2)
        $$ 
        where we may assume that the representation (12.2) is irreducible so that 
        $bP = \bigcup_{i=1}^n F_i$ where $F_i = \{z\in\C^N\colon\ 
        \Re(<z|w_i>)=1\}\cap P\ \ (1\leq i\leq n)$\ are the 
        facets of $P$.
        
        Since $P$ is compact there is an $R<\infty$ such that
        $$
        |<z|w_i>| \leq R\ \ (z\in P,\ \ 1\leq i\leq n).
        \eqno (12.3)
        $$
        By Proposition 12.1 there is an $\eta >0$ such that for each $j,\ 1\leq j\leq n$,
        $$
        bP\cap \{ z\in\C^N\colon \ 1-\eta<\Re (<z|w_j>) < 1\} \subset \cU .
        $$
        Passing to a smaller $\eta $ if necessary we may assume that 
        $$
        K\subset \{ z\in\C^N\colon \Re (<z|w_j>) \leq 1-\eta \} \hbox{\ \ for each\ \ }j, \ 1\leq j, \leq n .
        \eqno (12.4)
        $$
        By the remark following Proposition 12.1, for each $j,\ 1\leq j\leq n$ we have
        $$
        \bigcup_{i=1, i\not= j}^n [F_i\setminus\cU] \subset  \{ z\in\C^N\colon \Re (<z|w_j>) \leq 1-\eta \} .
        \eqno (12.5)
        $$
        Let $\varepsilon >0$ and $L<\infty $. 
        By the Runge theorem there is a polynomial 
        $\Phi\colon\ \C\rightarrow\C$  such that
        $$
        |\Phi  (\z)- (L+\varepsilon)|< \varepsilon/n\ \ (\z\in R\overline\D,\ \Re (\z)\geq 1)
        \eqno (12.6)
        $$
        $$|\Phi (\z )|<\varepsilon/n\ \ (\z\in R\overline\D,\ \Re (\z)\leq 1-\eta)
        \eqno (12.7)
        $$
        For each $j,\ 1\leq j\leq n$, consider the polynomial $f_j(z) =\Phi (<z|w_j>)$. By (12.4),
        $$
        |f_j(z)| <\varepsilon/n\ \ \ (z\in K)
        \eqno (12.8)
        $$
        and by (12.5) and (12.7), , 
        $$
        |f_j(z)|<\varepsilon/n\ \ \ \bigl( z\in\bigcup_{i=1, i\not=j}^n F_i\setminus\cU\bigr) .
        \eqno (12.9)
        $$
        Further, if $z\in F_j$ then $\Re(<z|w_j>)=1$  so by (12.6)
        $$
        |f_j(z)-(L+\varepsilon)|<\varepsilon/n\ \ (z\in F_j).
        \eqno (12.10)
        $$
        Now, let $f=\sum_{j=1}^nf_j$. If $1\leq j\leq n$ and if $z\in F_j\setminus\cU$ then by 
        (12.9) and (12.10)
        $|f(z)-(L+\varepsilon| \leq |f_j(z)-(L+\varepsilon)|+\vert 
        \sum_{i=1, i\not=j}^n f_i(z)\vert 
        \leq \varepsilon /n +(n+1)\varepsilon /n = \varepsilon$ 
        which implies that $\Re (f(z)\geq L \ (z\in F_j\setminus\cU, 1\leq j\leq n)$ 
        so $\Re (f(z))\geq L\ \ (z\in bP\setminus \cU)$. 
        Finally, by (2.8), \ $|f(z)|<\varepsilon \ \ (z\in K)$. This completes the proof.
        \vskip 2mm
        \bf\noindent Proof of Theorem 1.1\ \ \rm Let $P_n$ be the sequence of convex polytopes 
        from Theorem 2.1 and let
        $\theta_n$ be a decreasing sequence of positive numbers such that 
        $\sum_{n=1}^\infty \theta_n <\infty$. For each $n$, let $\cU_n\subset bP_n$ be the
        $\theta_n$-neighbourhood of $\hbox{skel}(P_n)$ in $bP_n$. The theorem will be proved once we 
        have constructed a holomorphic function $f$ on $\B_N$ such that
        $$
        \Re (f(z)) \geq n\ \ (z\in bP_n\setminus \cU_n,\ n\in\N) .
        \eqno (12.11)
        $$
        To see this, let $f$ satisfy (12.11) and suppose that $p\colon\ [0,1)\rightarrow \B_N$ is 
        a path such that $\lim_{t\rightarrow 1}|p(t)|=1$. Suppose that $f$ is bounded on $p([0,1))$. By 
        (12.11) there is some $n_0$ such that for each $n\geq n_0,\ \ p([0,1)) $ meets $bP_n$ only at 
        $\cU_n$. By
        Corollary 2.2 it follows that $p$ has infinite length.
        
        We shall construct a sequence $f_n$ of polynomials from $\C^N$ to $\C$ such
        that for each $n\in\N$, 
        $$
        \eqalign{ &(i)\ \Re (f_n(z)) \geq n+1 \ \  \hbox{on} \ \ bP_n\setminus \cU_n \cr
        &(ii)\ |f_{n+1}(z)-f_n(z)| \leq 1/2^{n+1} \hbox{\ on\ }P_n }.
        $$
        Suppose that we have done this. By (ii) the sequence converges uniformly on compacta in 
        $\B_N$ so the limit $f$ is 
        holomorphic on $\B_N$. If $z\in bP_n\setminus \cU_n$ then we have 
        $$
        f(z) = f_n(z)+\sum_{j=n}^\infty [f_{j+1}(z)-f_j(z)]
        $$
        so by (ii), \ \ $|f(z)-f_n(z)|<1$\ on\ $bP_n\setminus\cU_n$ and therefore 
        $\Re (f(z))\geq \Re (f_n(z)) -1\geq n$ on $bP_n\setminus \cU_n$ so that $f$ satisfies (12.11). 
        
        We construct $f_n$ by induction. Suppose that for some $m\in\N$ we have constructed 
        $f_m$  which satisfies
        $$
        \Re(f_m(z)) \geq m+1\ \ \hbox{on\ \ } bP_m\setminus \cU_m.
        $$
        Choose $T<\infty $ so large that
        $$\Re (f_m(z)) + T\geq m+2\ \ \hbox{\ on\ \ \ }bP_{m+1} .
        \eqno (12.12)
        $$
        By Lemma 12.2 there is a polynomial $g$ such that
        $$
        \Re(g(z))\geq T\ \ \hbox{on\ \ } bP_{m+1}\setminus\cU_{m+1}
        \eqno (12.13)
        $$
        and
        $$
        |g(z)| \leq (1/2)^{m+1}\ \ \hbox{on\ \ }P_m .
        \eqno (12.14)
        $$
        Put $f_{m+1}=f_m+g$. By (12.13) we have 
        $$
        \Re(f_{m+1}) = \Re(f_m+g) = \Re (f_m) + \Re (g) 
        \geq \Re (f_m)+T\geq m+2\hbox{\ on\ }bP_{m+1}\setminus \cU_{m+1}.
        $$
        and by (12.14) we have  $|f_{m+1}-f_m|<(1/2)^{m+1}\hbox{\ \ \ on\ \ }P_m$. 
        Theorem 1.1 is proved.
        \vskip 4mm
        \bf 13.\ Concluding remarks \rm
        \vskip 2mm
        We have proved Theorem 2.1 in $\R^M$ with $M\geq 3$. Theorem 2.1 holds 
        also in $\R^2$ where the proof is much simpler. 
        One can use a sequence of pairs of regular polygons.
        
        Having in mind the length of the proof of Theorem 2.1 one could say that 
        the principal result of the present paper is Theorem 2.1. It belongs to convex geometry
        and is not related to complex analysis. In its complex analysis consequence, 
        Theorem 1.1, the \it real part \rm 
        of the holomorphic function $f$ is unbounded on every path of finite length in $\B_N$ that ends 
        on $b\B_N$. Notice that by the maximum principle the zero sets of (real) pluriharmonic functions 
        on $\B_N,\ \ N\geq 2$, have no compact components. Applying Sard's theorem to the real part of the
        function $f$ obtained in Theorem 1.1 we get
        \vskip 2mm 
        \noindent \bf THEOREM 13.1\ \ \it Given $N\geq 2 $ there is a complete, 
        closed, real hypersurface of $\B_N$ which is the zero set of a (real) 
        pluriharmonic function on $\B_N$. \rm
        \vskip 2mm
        In the special case when $k=1$ and $N=2$ our Corollary 1.2 provides the 
        existence of a complete properly embedded 
        complex curve in $\B_2$. The existence of such a curve also follows from a recent paper
        of A.\ Alarc\'{o}n and J.\ F.\ L\'{o}pez [AL2]. Their 
        proof is completely different from the one presented here. However, neither of the proofs 
        provides any information about the topology of the curve so the following question remains open:
        \vskip 2mm
        \noindent \bf QUESTION 13.1\ \ \it Does there exist a complete proper holomorphic embedding 
        $f\colon\ \D\rightarrow \B_2$? \rm 
        \vskip 2mm
        \noindent Knowing now that for each $N\geq 2$ there are complete closed complex 
        hypersurfaces in $\B^N$ one may ask also
        \vskip 2mm
        \noindent \bf QUESTION 13.2\ \ \it Given $N\geq 2$, does there exist a complete 
        proper holomorphic 
        embedding $f\colon\ \B_N\rightarrow \B_{N+1}$ ?
        \vskip 10mm
        \noindent \bf ACKNOWLEDGEMENT\ \ \rm The author is grateful to David Eppstein 
        and John M.\ Sullivan for helpful 
        suggestions. He is also grateful to Toma\v z Pisanski for his interest.
        
        This work was supported by the Research Program P1-0291 from ARRS, Republic of Slovenia.
        \vfill
        \eject
        \centerline{\bf REFERENCES}
        \vskip 4mm
        \rm
        \noindent [AL1]\ A.\ Alarc\'{o}n and F.\ J.\ L\'{o}pez:\ Null curves in $\C^3$ and Calabi-Yau conjectures.
        
        \noindent Math\ Ann.\ 355 (2013) 429-455
        \vskip 1mm
        \noindent [AL2]\ A.\ Alarc\'{o}n and F.\ J.\ L\'{o}pez:\ Complete bounded complex curves in $\C^2$.
        
        \noindent Preprint, Arxiv:1305.2118v2 , May 2013
        \vskip 1mm
        \noindent 
        [AF]\ A.\ Alarc\'{o}n and F.\ Forstneri\v c:\ Every bordered Riemann surface is a complete proper curve in a ball.
        
        \noindent Math.\ Ann.\ 357 (2013) 1049-1070
        \vskip 1mm
        \noindent [B]\ A.\ Bronsted:\ \it An Introduction to Convex Polytopes.\rm
        
        \noindent Graduate Texts in Math.\ 90, 1983.  Springer-Verlag New York Inc.
        \vskip 1mm
        \noindent [CS]\ J.\ H.\ Conway, N.\ J.\ A.\ Sloane:\ \it Sphere Packings, Lattices 
        and Groups. \rm
        
        \noindent Grundl.\ Math.\ Wiss. 290, 1988 Springer Verlag New York Inc.
        \vskip 1mm
        \noindent [GS]\ J.\ Globevnik and E.\ L.\ Stout:\ Holomorphic functions with highly 
        noncontinuable boundary behavior.
        
        \noindent J.\ Anal.\ Math.\ 41 (1982) 211-216
        \vskip 1mm
        \noindent [J]\ P.\ W.\ Jones:\ A complete bounded complex submanifold of $\C^3$.
        
        \noindent Proc.\ Amer.\ Math.\ Soc.\ 76 (1979) 305-306
        \vskip 1mm
        \noindent [MUY] F.\ Martin,\ M.\ Umehara and K.\ Yamada:\ Complex bounded holomorphic curves 
        immersed in $\C^2$ with arbitrary genus.
        
        \noindent Proc.\ Amer.\ Math.\ Soc.\ 137 ( 2009) 3437-3450
        \vskip 1mm
        \noindent [R]\ W.\ Rudin:\ \it Function Theory in the Unit Ball of $\C^n$ \rm
        
        \noindent Grundl.\ Math.\ Wiss.\ 241, 1980. Springer-Verlag New York Inc.
        \vskip 1mm
        \noindent [Y1]\ P.\ Yang:\ Curvature of complex submanifolds of $\C^n$.
        
        \noindent J.\ Diff.\ Geom.\ 12 (1977) 499-511
        \vskip 1mm
        \noindent [Y2]\ P.\ Yang:\ Curvature of complex submanifolds of $\C^n$. 
        
        \noindent In: Proc.\ Symp.\ Pure.\ Math.\ Vol.\ 30, part 2, pp. 135-137. Amer.\ Math.\ Soc., 
        Providence, R.\ I.\ 1977
        \vskip 10mm
        \noindent Institute of Mathematics, Physics and Mechanics
        
        \noindent Ljubljana, Slovenia
        
        \noindent josip.globevnik@fmf.uni-lj.si
        
        \end